\documentstyle[12pt,leqno]{article}

\def\np{\noindent}
\def\bs{\bigskip}

\def\semi{\hbox{ $\times $ \kern-.972em \raise.12719em\hbox{ $_{^|}$}  }}

\def\be{\begin{enumerate}}
\def\ee{\end{enumerate}}
\def\bi{\begin{itemize}}
\def\ei{\end{itemize}}

\def\cR{{\cal R}}
\def\cV{{\cal V}}

\begin{document}

\centerline{\large{Addendum to ``Groups of Ribbon Knots''}}\bigskip

\centerline{ by Ka Yi Ng}

\

The purpose of this document is to clarify the inductive step described in the proof of 
Theorem 3.2 in \cite{ng}. In the second last sentence of the proof, it says, `it follows from the
inductive proof that $\cR_n$  is of index two.' The question of how this assertion is verified was first
raised by Dror Bar Natan and his student Ofer Ron \cite{barnatan}. \ To avoid any confusions that may arise
in the future, the author of the paper \cite{ng} would like to fill in the details to show that $\cR_n$  is
of index two. Readers of this document are assumed to have familiarity with \cite{ng}. For your reference,
Theorem 3.2 is stated below.\bs

\np {\bf  Theorem 3.2} {\it $\cR_n$    forms a subgroup of the free abelian group $\cV_n$  of index two. 
So its rank is the  same as the rank of $\cV_n$  and is the number of linearly independent primitive
rational invariants of order $\leq n$.} \bs

To show that $\cR_n$  is of index two, we shall establish the following two lemmas. \bs

\np {\bf Lemma A:} \ {\it  Let $D_1,\dots,D_r$  be all the distinct non-split chord diagrams of order $n\geq
2$.  Let $\lambda_1,\dots,\lambda_r$ be $r$  arbitrary  integers. Then for $n> 2$ , there exists a ribbon
knot
$K_n$  whose additive invariants of order $<n$ are trivial and $V(K_n) = \sum_{i=1}^r \lambda_i V(D_i)$ 
for all additive invariants $V$ of order $n$.  When $n=2$, there exists a ribbon knot $K_2$  with 
$V(K_2) = 2 \lambda_1 V(D_1)$  for all
additive invariants $V$  of order 2.}\bs

\np {\bf Lemma B:} \ {\it   Given a knot $K$  whose additive invariants of order $<n$  are trivial, there
exist n-chord  diagrams $D_1,\dots, D_r$  and integers $\lambda_1,\dots,\lambda_r$  such that 
$V(K) = \sum_{i=1}^r \lambda_i V(D_i)$   for all additive invariants \ $V$  of order $n$.}\bs

Lemma A is an extended version of Lemma 3.1 in \cite{ng} and can be proved using the same arguments as 
given in \cite{ng}. Lemma B is proved in \cite{ng-stan}. We shall now make use of these two lemmas 
to prove Lemma C. \bs

\np {\bf Lemma C:} {\it  Let $X$  be the right-handed trefoil. Let $K$  be an arbitrary knot and let
$n\geq 2$  be an integer.  Define $s$  to be 1 if the Arf invariant of the knot $K$  is non-zero and 0
otherwise. Then one can  construct a
ribbon knot $K_n$ so that $V(K \# sX \# K_n) = 0$  for all additive invariants $V$  of order
$\leq n$.}
\bs

\np {\bf  Proof of Lemma C:}  For $n=2$ , the knot can always be represented by an integral multiple
$\lambda$  of the non-split 2-chord  diagram $D$  so that $V(K) = \lambda V(D)$  for all additive
invariants of order 2. Write $\lambda$  as $2p-s$  where $p$  is an integer. By Lemma A, we can find a
ribbon knot $K_2$  with $V(K_2) = -2pV(D)$   for all additive invariants of order 2. The right-handed
trefoil knot has $V(X) = V(D)$  and thus, the case $n=2$  is proved.

Suppose a ribbon knot $K_n$  is constructed so that $V(K\# sX \# K_n) = 0$  for all additive
invariants of order $\leq n$.  Lemma B tells us that there exist integers $\lambda_1,\dots,\lambda_r$  and
$(n+1)$-chord diagrams $D_1,\dots ,D_r$  so that $V(K\# sX \# K_n) = \sum_{i=1}^r \lambda_i V(D_i)$ 
for all additive invariants of order $n+1$ . Apply Lemma A to construct a ribbon knot $R$  whose additive
inariants of order $\leq n$  are trivial and  $V(R) = -\sum_{i=1}^r \lambda_i V(D_i)$  for all additive
invariants of order $n+1$. Then  $V(K\# sX \# K_n \# R) = 0$ for all additive invariants of
order $\leq n+1$.

Lemma C implies that $V(K \# sX \# K_n) = 0$  for all rational-valued invariants of order
$\leq n$. According to the proof of  Theorem 3.2 in \cite{ng}, $\cR_n$  is a subgroup of ${\cal V}_n$.
Thus, as an element in
${\cal V}_n$, we have $[K]_n\in -[sX]_n + \cR_n$. This completes the proof that $\cR_n$  is of index two.

\

\noindent Ka Yi Ng, e-mail: Kayi.Ng@WallStreetSystems.com

\end{document}